\documentclass{article}
\usepackage[utf8]{inputenc}
\usepackage[english]{babel}
\usepackage{amsmath,amsthm,amsfonts}
\usepackage{amssymb}
\usepackage{graphicx} 
\usepackage{bm}
\usepackage{tikz}
\usepackage{multirow}
\usepackage{color}
\usepackage{placeins}
\usepackage{graphicx}
\usepackage{cite}

\usepackage{multicol}
\usepackage{wrapfig,lipsum,booktabs}

\graphicspath{ {./figs/} }
\newcommand{\revs}[1]{\textcolor{black}{#1}}

\begin{document}

\title{An Online Generalized Multiscale finite element method for heat and mass transfer problem with artificial ground freezing}% \thanks{Work is supported by the mega-grant of the Russian Federation Government (No. 14.Y26.31.0013 and  RSC No. 17–71–20055). }}

%\titlerunning{An Online Generalized Multiscale finite element method for unsaturated filtration problem in fractured media}
% If the paper title is too long for the running head, you can set
% an abbreviated paper title here
%

\author{
Denis Spiridonov
\thanks{Laboratory of Computational Technologies for Modeling Multiphysical and Multiscale Permafrost Processes, North-Eastern Federal University, 677000 Yakutsk,Republic of Sakha (Yakutia), Russia. 
Email: {\tt d.stalnov@mail.ru}.}
\and
Sergei Stepanov  \thanks{Laboratory of Computational Technologies for Modeling Multiphysical and Multiscale Permafrost Processes, North-Eastern Federal University, 677000 Yakutsk,Republic of Sakha (Yakutia), Russia
Email: {\tt cepe2a@inbox.ru}.}
\and
Vasiliy Vasil'ev \thanks{Department of Computational technologies, North-Eastern Federal University, 677000 Yakutsk,Republic of Sakha (Yakutia), Russia
Email: {\tt vasvasil@mail.ru}.}
}

%\and
%Maria Vasilyeva \thanks{Department of Mathematics and Statistics, Texas A\&M University, Corpus Christi, Texas, USA. 
%Email: {\tt maria.vasilyeva@tamucc.edu}.}

%\and
%Eric T. Chung 
%\thanks{Department of Mathematics,
%The Chinese University of Hong Kong (CUHK), Hong Kong SAR. 
%Email: {\tt tschung@math.cuhk.edu.hk}.}

%Yalchin Efendiev \thanks{Department of Mathematics \& Institute for Scientific Computation (ISC), Texas A\&M University, College Station, Texas, USA. Email: {\tt efendiev@math.tamu.edu}.} 
%\and
%Yalchin Efendiev \thanks{Department of Mathematics \& Institute for Scientific Computation (ISC),
%Texas A\&M University,
%College Station, Texas, USA. Email: {\tt efendiev@math.tamu.edu}.}

\maketitle
\maketitle              % typeset the header of the contribution
\begin{abstract}
The Online Generalized Multiscale Finite Element Method (Online GMsFEM) is presented in this study for heat and mass transfer problem in heterogeneous media with artificial ground freezing process. The mathematical model is based on the classical Stefan model, which depicts heat transfer with a phase change and includes filtration in a porous media. The model is described by a set of temperature and pressure equations. We employ a finite element method with the fictitious domain method to solve the problem on a fine grid. We apply a model reduction approach based on Online GMsFEM to derive a solution on the coarse grid. We can use the online version of GMsFEM to take less offline multiscale basis functions. We use decoupled offline basis functions built with snapshot space and based on spectral problems in our method. This is the standard approach of basis construction. We calculate additional basis functions in the offline stage to account for artificial ground freezing pipes. We use online multiscale basis functions to get a more precise approximation of phase change. We create an online basis that reduces error using local residual values. The accuracy of standard GMsFEM is greatly improved by using an online approach. Numerical results in a two-dimensional domain with layered heterogeneity are presented. To test the method's accuracy, we show results from a variety of offline and online basis functions. The results suggest that Online GMsFEM can deliver high-accuracy solutions with minimal processing resources.

\end{abstract}

\section{Introduction}
	
In this research, we investigate at a heat and mass transfer process with artificial ground freezing. Almost every aspect of building and mining uses artificially frozen filter soils\cite{andersland2003frozen, harris1995ground, newman2011artificial, alzoubi2017intermittent, jessberger1980theory}. For example, such technology is commonly utilized in mine sinking, tunneling, subway construction, and building construction, etc. Artificial freezing with cooling equipment is employed near piles in the construction of buildings on permafrost soils to assure stability by generating big chunks of frozen earth surrounding the pile, which will preserve the soil from defrosting throughout the summer \cite{loveridge2012energy, fei2019ground, vasilyeva2017reduced}. This practice has been demonstrated to strengthen building foundations in the Far North.

The authors describe mathematical models of heat and mass transfer processes in frozen and thawed soils, as well as artificial freezing of filter soils in the works \cite{esch2004thermal, alexiades2018mathematical}. A set of heat and mass transfer equations with a dynamic phase transition boundary is solved, which comprises a parabolic equation for temperature and an elliptic equation for pressure.  We employ a finite element method to approximate problem on the fine grid. \cite{bathe2007finite, szabo2021finite}. We apply the through counting approach to model temperature phase change \cite{burago2021numerical, rubinshteuin1971stefan} and the fictional domain method to model pressure phase change \cite{glowinski1994fictitious, glowinski1994fictitious1}.

These processes are generally modeled in very huge domains where the construction body is located.
 It is important to construct a highly fine computational grid in order to precisely describe the process at each point in a large region. This method demands high computational prices and high-end computing hardware, such as a supercomputer. As a way out of the situation, You can utilize homogenization methods \cite{talonov2016numerical, grigoriev2019effective, gavrilieva2018numerical, vasilyeva2017reduced} just design coarser grids to get out of the problem, but the accuracy of the solution will suffer. We propose a model reduction technique based on multiscale finite element method \cite{efendiev2009multiscale}. 

Multiscale finite element method(MsFEM) is particularly suited for modeling problems in highly heterogeneous regions \cite{hou1997multiscale, akkutlu2017multiscale}. Multiscale methods come in many forms, such as multiscale finite volume method (MsFVM)\cite{hajibeygi2008iterative, lunati2006multiscale, lunati2008multiscale, tyrylgin2019embedded} which uses the finite volume method to generate the basis functions \cite{eymard2000finite, moukalled2016finite}. Another MsFEM modification is the generalized finite element method(GMsFEM) \cite{chung2016adaptive, chung2014generalized, efendiev2013generalized} which builds several bases in each local domain by solving local spectral problems.  In the constraint energy minimizing generalized multiscale finite element method(CEM-GMsFEM), the basis building can be provided in an oversampled domain\cite{chung2018constraint, cheung2018constraint, fu2020constraint}. A multiscale method that satisfies the properties of mass conservation is called mixed multiscale finite element method (Mixed MsFEM) \cite{chen2003mixed, chung2015mixed, aarnes2008mixed}. An online generalized multiscale finite element method (Online GMsFEM) \cite{akkutlu2016multiscale, chung2015residual, tyrylgin2021multiscale, spiridonov2021online} is particularly suited for nonlinear problems because it executes the procedure of enriching a multiscale space during the online stage of the method. A special type of multiscale basis functions based on constrained energy minimization problems are developed in \cite{chung2018nonfrac, vasilyeva2019nonlocal, vasilyeva2019constrained, vasilyeva2019upscaling} and well-known as nonlocal multicontinuum method(NLMC).

We previously created a GMsFEM algorithm with an additional basis function for artificial ground freezing \cite{vasilyeva2020multiscale}. We expand our technique and employ multiscale online basis functions to predict phase change in the heat and mass transfer problem with artificial ground freezing in this study. We perform a model reduction procedure which includes offline and online stage. We create offline basis functions using local spectral problems in the offline stage. To better forecast processes in mediums with high contrast, we use snapshot space in basis computation. The considered medium features layered heterogeneity with a large value jump between two layers. In the computational domain, our problem involves frozen pipes. We employ offline additional basis functions to account  the effect of frozen pipes in order to make an accurate simulation.
%\revs{The efficiency of these multiscale basis functions has been shown in our previous work mentioned above \cite{vasilyeva2020multiscale}.} 
In the online stage, we construct an online multiscale basis functions. They assist in the accurate approximation of a shifting phase change boundary. Local residuals are used to compute online basis functions. \revs{Using online bases allows us to take smaller amount of offline basis functions with better accuracy.}

The numerical results are shown in a two-dimensional heterogeneous domain. We consider a two type of boundary condition for pressure. In the first example we set flow from right to left and in the second example we set flow from top to bottom. We complete the validation process by showing the results for a variety of offline and online basis functions.

The paper organized as follows. In Section 2, we present a mathematical model for heat and mass transfer problem with numerical algorithm for phase change boundary. In Section 3, we present an approximation on the fine grid using finite element method. Next, in Section 4, we describe an algorithm of online generalized multiscale finite element method. Numerical results for two test cases are presented in Section 5.

\section{Problem formulation}

We consider the heat and mass transfer model with artificial ground freezing. Freezing pipes within the soils enable artificial ground freezing process.
We use the classical Stefan model to simulate heat transfer processes with phase change \cite{samarskiy2009computational, Rubinstein1967problem}.  In this approximation, we assume that the phase change occurs at a given phase change temperature $T^*$. Let $\Omega^+(t)$ is the domain of the liquid phase where the temperature exceeds the phase change temperature:

\begin{equation}
\Omega^+(t) = {x| x \in \Omega, T(x,t)>T^*}, \nonumber  
\end{equation}

and $\Omega^-(t)$ is the domain of frozen phase:
\begin{equation}
\Omega^-(t) = {x| x \in \Omega, T(x,t)<T^*}, \nonumber  
\end{equation}

We denote by $\rho^+$,$c^+$ and $\rho^-$,$c^-$ the density and specific heat capacity of the liquid and frozen zones, respectively. We define the indicator function as piecewise defined function
\revs{
\begin{equation}
\phi = \phi(T) = 
\begin{cases}
0, & T<T^*, \\
1, & T>T^*.
\end{cases} \nonumber
\end{equation}}

For the coefficients of heat capacity and thermal conductivity, we have
\begin{equation}
\alpha(\phi) = \rho^- c^- + \phi(\rho^+ c^+ - \rho^- c^-), \ \ \ k(\phi) = k^-+\phi(k^+ - k^-).\nonumber
\end{equation}

Thermal processes, which are accompanied by phase change, absorption and release of latent heat of fusion, are described by the equation
\begin{equation} \label{eq1}
\Big(\alpha(\phi)+\rho^+L\frac{\partial \phi}{\partial T}\Big)\Big(\frac{\partial T}{\partial t}+u \nabla T\Big)-\nabla\cdot (k(\phi)\nabla T)=Q.
\end{equation} 

Here $L$ is the specific heat of the phase change, $u$ is the filtration speed in the soil.

We have heat capacity coefficients because we are considering heat transfer in a porous material:
\begin{equation}
c^- \rho^- =(1-m)c_{sc} \rho_{sc} + mc_l \rho_l, \ \ \ c^+\rho^+=(1-m)c_{sc}\rho_{sc}+mc_w \rho_w, \nonumber 
\end{equation}

\revs{here $m$ is the porosity and} the subscripts $sc$, $w$ , and $l$ denote the skeleton of porous medium, water and ice, respectively. We use comparable terms for heat conductivity coefficients in the solid and frozen zones:
\revs{
\begin{equation}
\lambda^-=(1-m)\lambda_{sc}+m\lambda_l, \ \ \ \lambda^+=(1-m)\lambda_{sc} + m\lambda_{w}. \nonumber
\end{equation}}

To take into account the filtration in the soil, we write the continuity equation for a completely saturated porous medium
\begin{equation} \label{eq2}
\frac{\partial}{\partial t}(\rho _w m)+\nabla \cdot(\rho _w u)=0.
\end{equation}
We use Darcy's law to express the relationship between the filtration rate $u$ and the pressure gradient:
\begin{equation} \label{eq3}
u=-\frac{k}{\mu}(\nabla p +\rho _w g),
\end{equation}
where $k$ is the absolute permeability tensor of the porous medium, $\mu$ is the fluid viscosity, $g$ is the free fall acceleration \revs{and $\rho _w$ is the density of water}.

Substituting \eqref{eq3} into \eqref{eq2} and neglecting compressibility for simplicity, we obtain equations for determining reservoir pressure
\begin{equation}\label{eq4}
-\nabla \cdot (\lambda \nabla p) = F, \ \ \ x\in\Omega^+,
\end{equation}
where $\lambda = \rho _w k/\mu$.

The system of equations \eqref{eq1}, \eqref{eq4} is the base for modeling the processes of thermal stabilization of filter soils. In our problem the main difficulties of numerical simulation are generated by the phase change, we need to solve the filtration problem with a free (unknown) boundary. We use computational algorithms for through counting, which were used by many authors in solving similar problems of heat and mass transfer with phase transformations \cite{samarskii1993numerical, belhamadia2012enhanced, bhattacharya2002fixed}.

In our implementation we assume that the phase change occurs not when $T^*=T$, but in some small interval near the temperature of phase change $[T^*-\delta,T^*+\delta]$. Instead of the indicator function $\phi$, we take its piecewise linear approximation of the following form:
\begin{equation}
\phi_\delta(T)=\begin{cases}
0, & T\leq T^*-\triangle, \\
\frac{T-T^*+\triangle}{2\triangle}, & T^* - \triangle < T < T^*+\triangle, \\
1, & T \geq T^*+\triangle.
\end{cases} \nonumber
\end{equation}
Then, we have
\revs{
\begin{equation}
\frac{d\phi_\triangle}{dT}=\begin{cases}
0, & T < T^*-\triangle, \\
\frac{1}{2\triangle}, & T^* - \triangle < T < T^*+\triangle, \\
0, & T > T^*+\triangle.
\end{cases} \nonumber
\end{equation}}

Thus, instead of \eqref{eq1} in the computational domain $\Omega$, we solve the following equation for temperature:
\begin{equation}\label{eq8}
\Big(\alpha(\phi_\triangle)+\rho^+L\frac{d\phi_\triangle}{dT}\Big)\Big(\frac{\partial T}{\partial t}+u \nabla T\Big)-\nabla\cdot(k(\phi_\triangle)\nabla T)=Q
\end{equation}

This non-linear parabolic equation is supplemented with appropriate initial and boundary conditions.

Since our problem \eqref{eq4}-\eqref{eq8} has a moving phase transition boundary, we are faced with the problem of constructing a numerical algorithm for calculating the pressure. For the numerical solution of such a problem without rebuilding the computational grid, we use the method of fictitious domain, which is based on the transition to solving the problem in a wider domain. An approximate solution, depending on the continuation parameter $\varepsilon$ , we will solve in the entire computational domain $\Omega$. When using a variant of the method of fictitious regions with continuation by leading coefficients, the solution is determined from the equation
\begin{equation}\label{eq9}
\nabla\cdot(\lambda_\varepsilon\nabla p)=F, \ \ \ x\in\Omega.
\end{equation}
Here, the discontinuous coefficient $\lambda_\varepsilon(x)$ is determined by the expression
\begin{equation}
\lambda_\varepsilon(x)=\begin{cases}
\lambda, \ \ \ x\in \Omega^+, \\
\lambda\varepsilon, \ \ \ x\in \Omega^- 
\end{cases}
\end{equation}
with very small $\varepsilon$. With this choice of coefficients, equation \eqref{eq9} simulates filtration in the region $\Omega^-$ with a very small coefficient $\lambda_\varepsilon(x) \rightarrow 0$ with the parameter $\varepsilon \rightarrow 0$.

We supplement our system \eqref{eq8},\eqref{eq9} with appropriate initial and boundary conditions. For temperature, we apply the initial condition:
\begin{equation} \label{eq5}
T=T_0(x), \ \ \ x\in\Omega.
\end{equation}

As the boundary condition for the temperature, we take the zero Neumann boundary condition over the entire boundary of the computational domain:

\begin{equation} \label{eq6}
-k\frac{\partial T}{\partial n}=0, \ \ \ x\in\partial\Omega.
\end{equation}

In our problem we have an artificial ground freezing process. To simulate it, we set the following boundary condition on the pipes:
\begin{equation}\label{eq10}
T=T_p,  \ \ \ x \in \Gamma _p
\end{equation}
where $\Gamma_p$ is the boundary of artificial ground freezing pipes. 

For the pressure, we set boundary conditions on the boundary of the entire computational domain $\partial \Omega = \Gamma_N\cup\Gamma_D$:
\begin{equation} \label{eq7}
-\lambda_\varepsilon  \frac{\partial p}{\partial n}=0, \ \ \ x\in \Gamma_N, \ \ \  p=q_D, \ \ \ x\in\Gamma_D,
\end{equation}

When using the computational algorithm for through counting, it is necessary to pay special attention to setting the parameters $\triangle$ and $\varepsilon$. \revs{These parameters indicate the width of the phase transition, they determine the accuracy of the phase change boundary for solving on a fine grid. $\triangle$ and $\varepsilon$ depend on the dimension of the grid, for an accurate approximation, a sufficient number of grid elements within the phase transition is required. In the paper, we chose these parameters based on the works \cite{vabishevich2014numerical, vasil2020accurate}, where the authors studied the influence of hyperparameters on the solution of the problem with a phase change.}

\section{Fine grid approximation}

Let us discretize in space using the finite element method of the initial-boundary value problem \eqref{eq8},\eqref{eq9}, \eqref{eq5}-\eqref{eq7} for the resulting basic system of equations, which is used to simulate artificial freezing processes. To reformulate the original system of temperature and pressure equations, we multiply equation \eqref{eq8} by the test functions $v$, and equation \eqref{eq9} by the test function $w$. We integrate our system using Green's formula and boundary conditions (6), (7) and obtain the following formulation:
\begin{equation}\label{eq11}
\begin{split}
\int_\Omega \Big(\alpha(\phi_\triangle)+\rho^+L\frac{d\phi_\triangle}{dT}\Big)\Big(\frac{\partial T}{\partial t}+u \nabla T\Big)v dx + \int_\Omega (k(\phi_\triangle)\nabla T &, \nabla v)dx =\int_\Omega Q v dx, \ \ \forall v\in H^1_0(\Omega),\\
\int_\Omega (\lambda_\varepsilon \nabla p, \nabla w)dx = \int_\Omega F w dx, & \ \ \forall w\in H^1_0(\Omega). 
\end{split}
\end{equation}
Here $H^1(\Omega)$ is  Sobolev space that consist of functions q such that $q^2$ and \revs{$|\nabla q|^2$} have a finite integral in $\Omega$ and $H^1_0(\Omega)={q\in H^1(\Omega): q(x)=0, \ x\in\Gamma_P}$.

To approximate the equation for temperature in time, we use the standard implicit scheme. The simplest linearization is used when specifying the coefficients from the previous time layer.\revs{Let $t_{max}=n\tau, \ n=1,2,...,N_t$ and $\tau$ be the time step.} We denote the solution at time $t_n$ by $T^n, \ p^n$. According to \eqref{eq11} we can derive a variational formulation:
\begin{equation}\label{eq12}
\begin{split}
\int_\Omega \Big(\alpha(\phi^n_\triangle)+\rho^+L\frac{d\phi_\triangle}{dT}(T^n)\Big) & \Big(\frac{T^{n+1}-T^n}{\tau}+(\lambda^n_\varepsilon \nabla p^n \nabla T^n)\Big)v dx +\\
+ \int_\Omega &(k(\phi^n_\triangle)\nabla T^{n+1} , \nabla v)dx =\int_\Omega Q v dx, \ \ \forall v \in H^1_0(\Omega),\\
\int_\Omega &(\lambda^n_\varepsilon\nabla p^{n+1}, \nabla w)dx = \int_\Omega F w dx, \ \ \forall w\in H^1_0(\Omega) 
\end{split}
\end{equation}

Let write approximation in the matrix form:

\revs{
\begin{equation}\label{eq122}
\begin{split}
S^n_T \frac{T^{n+1}-T^n}{\tau}+A^n_T T^{n+1} = G_T, \\
A^n_p p^{n+1} = G_p.
\end{split}
\end{equation}} 

Approximation in space \eqref{eq12} is associated with the choice of the corresponding finite element spaces for temperature and pressure and is carried out in the standard way. For both temperature and pressure, we use Lagrangian finite elements of the first order on a triangular mesh.

\section{Coarse grid approximation}

To make approximation on the coarse grid, we apply an Online Generalized Multiscale Finite Element method. A coarse grid $mathcalT H$ is denoted by
\begin{equation}
\mathcal{T}_H = \bigcup_j K_j, \nonumber
\end{equation}
here $K_j$ is the cell of coarse grid. Following that, we create a local coarse neighborhood domain $\omega _i$ it is produced by merging all coarse cells around one coarse grid vertex:
\begin{equation}
\omega _i=\bigcup_j{K_j \in \mathcal{T} \ : \ x_i \in K_j}, \nonumber
\end{equation}
where $i=\overline{1,N_c}$ and $N_c$ -- the number of coarse grid vertices.

We can identify two stages in GMsFEM: online and offline:

\ \ \ \ \textit{Offline stage:}

\begin{enumerate}
\item Creating local coarse neighborhood domains $\omega _i$ from the coarse grid;
\item In each local domain $\omega _i$ we construct an offline multiscale basis functions;
\item Construction an additional basis function for temperature $T$ in each local domain $\omega _i$ that contains ground freezing pipes;
\item Assembling a projection matrix $R$ of the multiscale space. This mathix consist of offline basis functions.
\end{enumerate}

\ \ \ \ \textit{Online stage:}

\begin{enumerate}
\item Assembling a system of equations on fine grid and projecting it on the multiscale space by $R$ matrix;
\item Solving a system of equations on the coarse grid by the resulting approximation on a coarse grid;
\item Constructing an online multiscale basis functions on certain time layers;
\item Resolving system on the coarse grid using new matrix $R^{on}$ with online bases.
\end{enumerate}
We should note, that we perform steps 3 and 4 only on certain time steps, otherwise we skip this steps. We add online bases not on every time layer, for example, we can do it on every 5-th or 10-th time layers. Following that, we go over the GMSFEM algorithm in further detail.

\subsection{Offline stage}

\textbf{Snapshot space.}

First, to build an offline multiscale basis functions, we need to define a snapshot space for temperature $V^{snap}_T$ and for pressure $V^{snap}_p$. Therefore, we have to build projection matrices for these spaces. We solve the following problems in each local domain $\omega _i$ to identify the matrix's members:

\begin{itemize}
\item For temperature:
\begin{equation} \label{eq14} 
\begin{cases}
-\nabla \cdot (k^+ \nabla \varphi^i_{T,j})=0, \ \ x\in \omega _i, \\
\varphi^i_{T,j} = \delta_j, \ \ x\in \partial \omega _i.
\end{cases}
\end{equation}
\item For pressure:
\begin{equation} \label{eq15} 
\begin{cases}
-\nabla \cdot (\lambda \nabla \varphi^i_{p,j})=0, \ \ x\in \omega _i, \\
\varphi^i_{p,j} = \delta_j, \ \ x\in \partial \omega _i.
\end{cases}
\end{equation}
where $\delta _j$ is the discrete delta function which equal 1 at the $j$-th fine grid node $x = x_j$ and zero elsewhere $(j = 1, ..., J_i$, $J_i$ is a number of fine grid vertices on boundary $\partial \omega _i$). And $\partial \omega _i$ is the outer boundary of local domain $\omega _i$.

The purpose of creating a snapshot space must be discussed. Almost all applied problems are modeled in domains with high coefficient contrast. Such domains may contain fractures or channels. We will compute our problem in domain with high coefficient contrast. We also need to employ a snapshot space to account all important properties of the solution and produce a good approximation. The snapshot vectors can describe most important characteristics of computational domain and help to construct more accurate projection to multiscale space.
\end{itemize}

Now we can define a snapshot spaces and projection matrices for temperature and pressure:
\begin{equation}\label{eq16}
\begin{split}
V^{snap,i}_T = \mbox{span}\{\varphi^i_{T,1}, ... , \varphi^i_{T,J_i}\}, \ \ \  R^{snap,i}_T=(\varphi^i_{T,1}, ... , \varphi^i_{T,J_i})^T, \\
V^{snap,i}_p = \mbox{span}\{\varphi^i_{p,1}, ... , \varphi^i_{p,J_i}\}, \ \ \  R^{snap,i}_p=(\varphi^i_{p,1}, ... , \varphi^i_{p,J_i})^T.
\end{split}
\end{equation}

\textbf{Spectral problem.}

The next step of offline multiscale basis construction is computing spectral problems in each local domain $\omega _i$. We should mention that we build a decoupled multiscale basis functions. To obtain a set of basis functions, we solve next eigenvalue problems:
\begin{itemize}
\item For temperature:
\begin{equation}\label{eq17}
\tilde{A}^i_T \tilde{\Psi^i}_{T,j} = \tilde{\lambda}\tilde{S}^i_T \tilde{\Psi^i}_{T,j},
\end{equation} 
with
\begin{equation}
\tilde{A}^i_T=R^{snap,i}_T A^i_T (R^{snap,i}_T)^T, \ \ \ \tilde{S}^i_T=R^{snap,i}_T S^i_T (R^{snap,i}_T)^T, \nonumber
\end{equation}
where
\begin{equation}
\begin{split}
A^i_T ={a_{ln}}, \ \ \ a_{ln}=\int_{\omega _i} (k^+ \nabla \psi_{m,l}, \nabla \psi_{m,n}) dx \\
S^i_T ={a_{ln}}, \ \ \ a_{ln}=\int_{\omega _i} (k^+ \psi_{m,l}, \psi_{m,n}) dx \\
\end{split} \nonumber
\end{equation}
\item For pressure
\begin{equation}\label{eq18}
\tilde{A}^i_p \tilde{\Psi^i}_{p,j} = \tilde{\lambda}\tilde{S}^i_p \tilde{\Psi}^i_{p,j},
\end{equation} 
with
\begin{equation}
\tilde{A}^i_p=R^{snap,i}_p A^i_p (R^{snap,i}_p)^T, \ \ \ \tilde{S}^i_p=R^{snap,i}_p S^i_p (R^{snap,i}_p)^T, \nonumber
\end{equation}
where
\begin{equation}
\begin{split}
A^i_p ={a_{ln}}, \ \ \ a_{ln}=\int_{\omega _i} (\lambda \nabla \psi_{m,l}, \nabla \psi_{m,n}) dx \\
S^i_p ={a_{ln}}, \ \ \ a_{ln}=\int_{\omega _i} (\lambda \psi_{m,l}, \psi_{m,n}) dx \\
\end{split} \nonumber
\end{equation}
\end{itemize}
Then we obtain solutions of spectral problems on $\omega _i$ space $\Psi^i_{T,j}=R^{snap,i}_T \tilde{\Psi}^i_{T,j}$, $\Psi^i_{p,j}=R^{snap,i}_p \tilde{\Psi}^i_{p,j}$. We generate offline multiscale basis functions using the initial $M i$ eigenvectors that correspond to the smallest $M i$ eigenvalues.

\textbf{Multiscale space.}

\begin{figure}[h!]
\begin{center}
\includegraphics[width=0.8\linewidth]{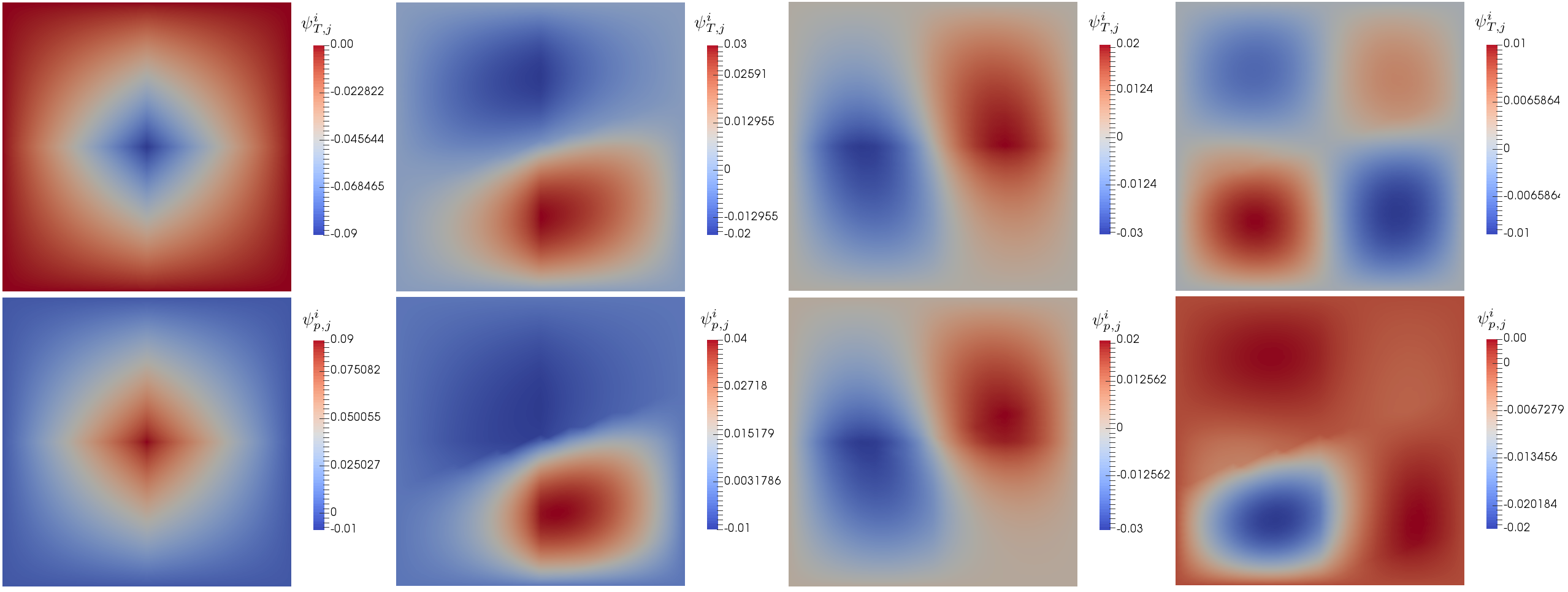}
\end{center}
\caption{The example of 4 offline multiscale basis functions in $\omega _i$. First row: bases for temperature. Second row: bases for pressure.}
\label{offbasis}
\end{figure}

To get offline bases, we utilize the following formula:
\begin{equation}\label{eq20}
\psi^i_{T,j}=\Psi^i_{T,j} \chi^i, \ \ \ \psi^i_{p,j}=\Psi^i_{p,j} \chi^i.
\end{equation}
Where $\chi^i$ is the linear partition of unity function that equals to standard coarse grid nodal basis function and $j=1,...,M_i$. The example of first 4 offline multiscale basis functions in one local domain $\omega _i$ is presented in Figure \ref{offbasis}.

After that, we'll be able to create multiscale spaces
\begin{equation}
\begin{split}
V^T_{ms}=\mbox{span}\{\psi^1_{T,1},...,\psi^1_{T,N_c},...,\psi^{N_c}_{T,1},...,\psi^{N_c}_{T,N_c}\}, \\
V^p_{ms}=\mbox{span}\{\psi^1_{p,1},...,\psi^1_{p,N_c},...,\psi^{N_c}_{p,1},...,\psi^{N_c}_{p,N_c}\},
\end{split} \nonumber
\end{equation}
and projection matrices
\begin{equation}
\begin{split}
R_T=(\psi^1_{T,1},...,\psi^1_{T,N_c},...,\psi^{N_c}_{T,1},...,\psi^{N_c}_{T,N_c})^T, \\
R_p=(\psi^1_{p,1},...,\psi^1_{p,N_c},...,\psi^{N_c}_{p,1},...,\psi^{N_c}_{p,N_c})^T.
\end{split} \nonumber
\end{equation}

\textbf{Additional offline multiscale basis function.}

\begin{figure}[h!]
\begin{center}
\includegraphics[width=0.8\linewidth]{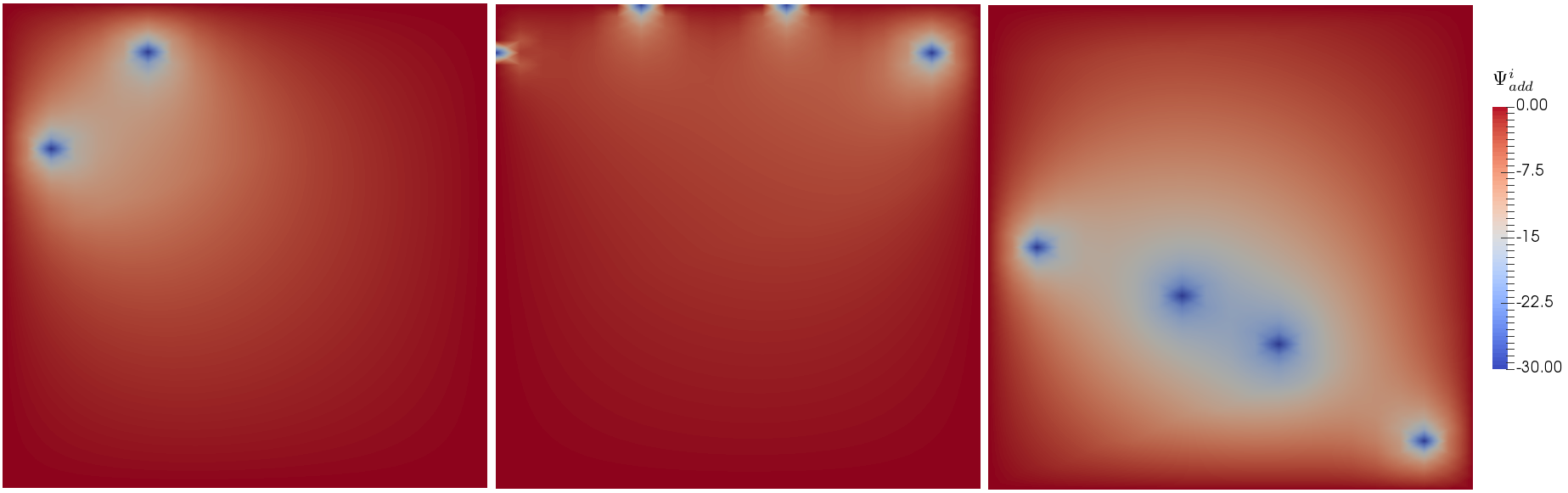}
\end{center}
\caption{The solution of local problem \ref{eq20} in three local domains $\omega _i$.}
\label{basisadd}
\end{figure}

In our problem, we have freezing  pipes in modeling of temperature field. We must consider the impact of boundary condition $\Gamma _P$. To calculate additional basis functions for temperature, we solve the next local problem in each local domain that includes freezing pipes
\begin{equation}\label{eq20}
\begin{cases}
-\nabla \cdot (k^+ \nabla \Psi^i_{add})=0, \ \ x\in \omega _i, \\
\Psi^i_{add} = 0, \ \ x\in \partial \omega _i, \\
\Psi^i_{add} = T_p, \ \ x\in \Gamma _P. 
\end{cases}
\end{equation} 
Then we multiply the solution of local problem \eqref{eq20} on partition of unity function to obtain an additional basis function $\psi^i_{add}=\Psi^i_{add}\chi^i$. \revs{It should be noted, that the calculation of additional bases does not depend on the location of the freezing pipes, they can be both inside the coarse block and coincide with the coarse edges. One can see it in the solution of problem \eqref{eq20} in Figure \ref{basisadd}.}

As a result, we get the following multiscale space and projection matrix for temperature:
\begin{equation}
\begin{split}
V^T_{ms}=\mbox{span}\{\psi^1_{T,1},...,\psi^1_{T,N_c},...,\psi^{N_c}_{T,1},...,\psi^{N_c}_{T,N_c}, \psi^1_{add},..,\psi^{N_p}_{add}\}, \\
R_T=(\psi^1_{T,1},...,\psi^1_{T,N_c},...,\psi^{N_c}_{T,1},...,\psi^{N_c}_{T,N_c},\psi^1_{add},..,\psi^{N_p}_{add})^T, \\
\end{split} \nonumber
\end{equation}
where $N_p$ is the number of $\omega_i$ with ground freezing pipes.

\subsection{Online stage}

One can see that in offline basis construction we didn't take into account the phase change. We use thermal conductivity $k$ and porous medium permeability $\lambda$ from liquid zone. We can't take into account the phase change by offline bases, because it changes by time.  We can neglect the phase transition, as we did earlier in \cite{vasilyevastep2020multiscale}, but in this paper we will take into account the phase change by online multiscale basis functions. The Online GMsFEM is well suited for solving non-linear problems. This method can significantly reduce the number of offline basis functions with good improvement in accuracy. It is better to add 1 or 2 online basis in each $\omega _i$ then computing a large number of offline bases. But this process has a big computational cost, because of this we compute online bases only in certain time layer, for example, in each 5-th or 10-th time layer. We construct an online multiscale basis function based on local residuals that provides fast error decay.

To make update of projection matrices in each $n$-th time step, firstly, we need to solve the coarse grid system using only offline multiscale basis functions to find local residuals:

\revs{
\begin{equation}\label{eq22}
\begin{split}
S^n_{T,c} \frac{T^{n+1}_c-T^n_c}{\tau}+A^n_{T,c} T^{n+1}_c = G_{T,c}, \\
A^n_{p,c} p^{n+1}_c = G_{p,c},
\end{split}
\end{equation} }
where
\revs{
\begin{equation}
\begin{split}
S^n_{T,c}=R^n_T S^n_T (R^n_T)^T, \ \ A^n_{T,c}=R^n_T A^n_T (R^n_T)^T, \ \ G_{T,c} = R^n_T G_T, \ \ T^{n+1}_{ms} = (R^n_T)^T T^{n+1}_c, \\
 A^n_{p,c}=R^n_p A^n_p (R^n_p)^T, \ \ G_{p,c} = R^n_p G_p, \ \ p^{n+1}_{ms} = (R^n_p)^T p^{n+1}_c,
\end{split} \nonumber
\end{equation}}
where $T^{n+1}_{ms}$, $p^{n+1}_{ms}$ are multiscale solutions that reconstructed on fine grid, $R^n_T$ and $R^n_p$ are predefined projection matrices without enrichment with online bases at $n$-th time layer 
\begin{equation}
\begin{split}
R^n_T=R_T=(\psi^1_{T,1},...,\psi^1_{T,N_c},...,\psi^{N_c}_{T,1},...,\psi^{N_c}_{T,N_c},\psi^1_{add},..,\psi^{N_p}_{add})^T, \\
R^n_p=R_p=(\psi^1_{p,1},...,\psi^1_{p,M_{N_c}},...,\psi^{N_c}_{p,1},...,\psi^{N_c}_{p,M_{N_c}})^T.
\end{split} \nonumber
\end{equation}

To complete online stage with solving coarse system with updating online multiscale bases functions, firstly, we solve the coarse scale system \eqref{eq22} using offline projection matrices $R^n_T$ and $R^n_p$. Then we add online bases in these matrices and obtain new projection matrices \revs{$R^{(L,n)}_T$, $R^{(L,n)}_p$}. We use these matrices until the next update procedure, we update projection matrices on each certain $n$-th time layer. The process can be iterable to add more multiscale online basis functions:
\begin{equation}
\begin{split}
R^n_T=R^{L,n}_T=(\psi^1_{T,1},...,\psi^1_{T,N_c},...,\psi^{N_c}_{T,1},...,\psi^{N_c}_{T,N_c},\psi^1_{add},..,\psi^{N_p}_{add}, \vartheta^1_{T,l},...,\vartheta^{N_c}_{T,l},...,\vartheta^1_{T,L},...,\vartheta^{N_c}_{T,L})^T, \\
R^n_p=R^{L,n}_p=(\psi^1_{p,1},...,\psi^1_{p,M_{N_c}},...,\psi^{N_c}_{p,1},...,\psi^{N_c}_{p,M_{N_c}},\vartheta^1_{p,l},...,\vartheta^{N_c}_{p,l},...,\vartheta^1_{p,L},...,\vartheta^{N_c}_{p,L})^T,
\end{split} \nonumber
\end{equation}
where $L$ denotes the number of online iterations. On the time layer where we don't update online bases, we use $R^n_T$, $R^n_p$ with old online bases.

\begin{figure}[h!]
\begin{center}
\includegraphics[width=0.8\linewidth]{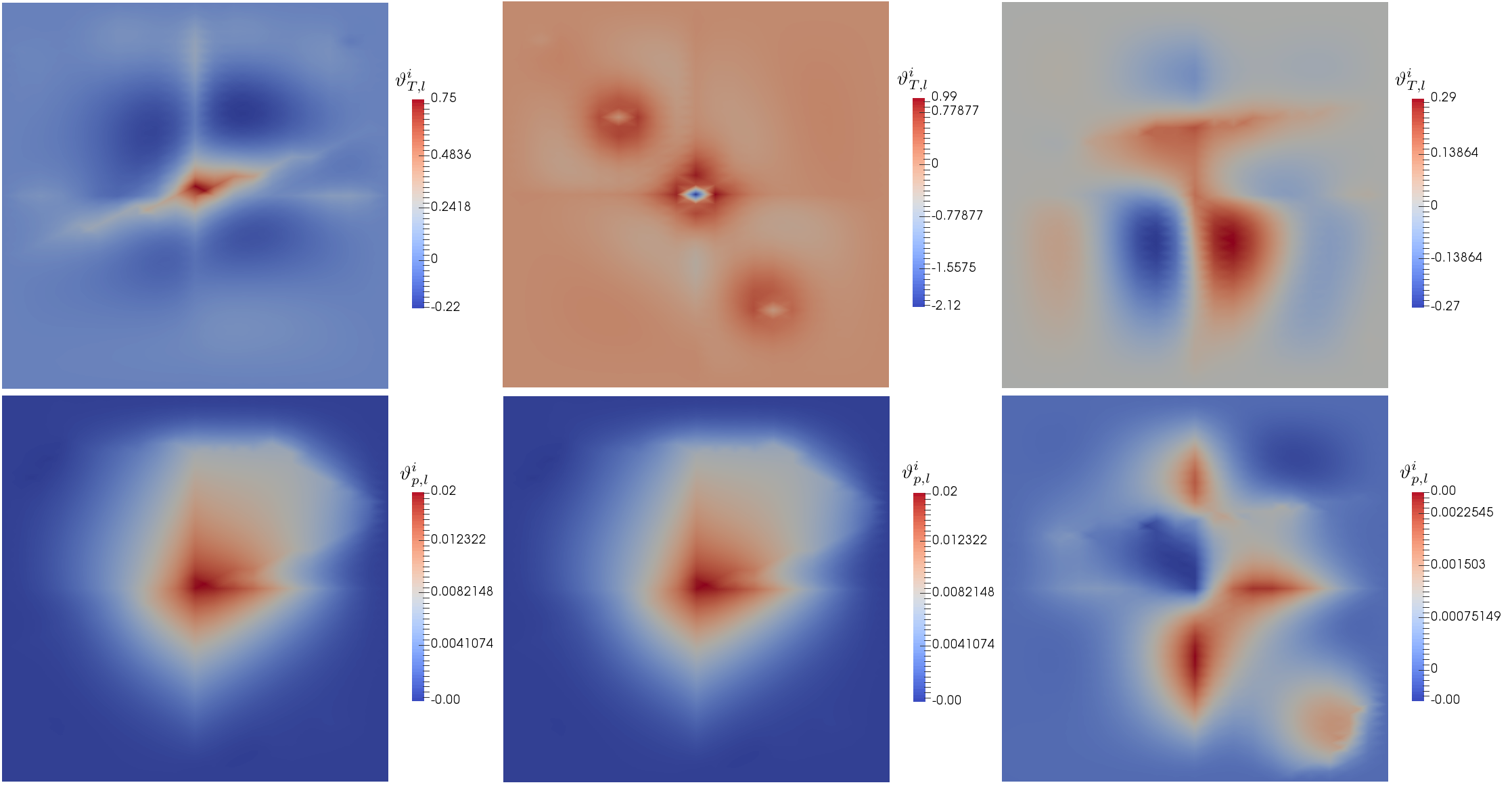}
\end{center}
\caption{An example of online multiscale basis functions in three local domains $\omega _i$. First row: bases for temperature. Second row: bases for pressure.}
\label{basisonline}
\end{figure}

Next we present a construction of online multiscale basis functions that based on solution of the following problem in each local domain $\omega _i$:
\begin{equation} \label{eq24}
\begin{split}
a_{T,\omega _i}(\Phi^i_{T,l},v) = r^l_{T,\omega _i}(v), \\
a_{p,\omega _i}(\Phi^i_{p,l},w) = r^l_{p,\omega _i}(w), 
\end{split} \ \ \ \ \ l=1,...,L,
\end{equation}
where
\begin{equation}
\begin{split}
a_{T,\omega _i}(\Phi^i_{T,l},v) = \int_\Omega \Big(\alpha(\phi^n_\triangle)+\rho^+L\frac{d\phi_\triangle}{dT}(T^n_{ms})\Big) & \Big(\frac{\Phi^i_{T,l}}{\tau}+(\lambda^n_\varepsilon \nabla p^n_{ms} \nabla T^n_{ms})\Big)v dx +\\
+ \int_\Omega &(k(\phi^n_\triangle)\nabla \Phi^i_{T,l} , \nabla v)dx, \\
r^l_{T,\omega _i}(v) = \int_\Omega Q v dx - \int_\Omega \Big(\alpha(\phi^n_\triangle)+\rho^+L\frac{d\phi_\triangle}{dT}(T^n_{ms})\Big) & \Big(\frac{T^{n+1}_{ms}-T^n_{ms}}{\tau}+(\lambda^n_\varepsilon \nabla p^n_{ms} \nabla T^n_{ms})\Big)v dx -\\
- \int_\Omega &(k(\phi^n_\triangle)\nabla T^{n+1}_{ms} , \nabla v)dx, \\
a_{p,\omega _i}(\Phi^i_{p,l},w) = \int_\Omega &(\lambda^n_\varepsilon\nabla \Phi^i_{p,l}, \nabla w)dx, \\
r^l_{p,\omega _i} = \int_\Omega F w dx - \int_\Omega &(\lambda^n_\varepsilon\nabla p^{n+1}_{ms}, \nabla w)dx,
\end{split} \nonumber
\end{equation}
with the boundary condition $\Phi^i_{T,l}=0, \Phi^i_{p,l}=0$ on $\partial \omega _i$. After solving problem \eqref{eq24} we multiply the solution on partition of unity functions to obtain an online multiscale basis functions:
\begin{equation}\label{eq25}
\vartheta^i_{T,l}=\Phi^i_{T,l} \chi^i, \ \ \ \vartheta^i_{p,l}=\Phi^i_{p,l} \chi^i.
\end{equation}
We demonstrate an example of online multiscale basis functions in Figure \ref{basisonline}.

Next we update multiscale spaces with online basis functions:
\begin{equation}
\begin{split}
V^T_{ms}=\mbox{span}\{ \psi^1_{T,1},...,\psi^1_{T,N_c},...,\psi^{N_c}_{T,1},...,\psi^{N_c}_{T,N_c},\psi^1_{add},..,\psi^{N_p}_{add}, \vartheta^1_{T,l},...,\vartheta^{N_c}_{T,l} \}, \\
V^p_{ms}=\mbox{span}\{ \psi^1_{p,1},...,\psi^1_{p,M_{N_c}},...,\psi^{N_c}_{p,1},...,\psi^{N_c}_{p,M_{N_c}},\vartheta^1_{p,l},...,\vartheta^{N_c}_{p,l} \},
\end{split} \ \ \ l=1,...,L. \nonumber
\end{equation}

Next, we use current solutions $T^{l,n}_{ms}$, $p^{l,n}_{ms}$ to find a certain number of online basis functions in each local domain $\omega_i$:

\revs{
\begin{equation}\label{eq26}
\begin{split}
S^{l,n}_{T,c} \frac{T^{l,n+1}_c-T^{l,n}_c}{\tau}+A^{l,n}_{T,c} T^{l,n+1}_c = G_{T,c}, \\
A^{l,n}_{p,c} p^{l,n+1}_c = G_{p,c},
\end{split}
\end{equation} 
}
where
\revs{
\begin{equation}
\begin{split}
S^{l,n}_{T,c}=R^{l,n}_T S^{l,n}_T (R^{l,n}_T)^T, \ \ A^{l,n}_{T,c}=R^{l,n}_T A^{l,n}_T (R^{l,n}_T)^T, \ \ G_{T,c} = R^{l,n}_T G_T, \ \ T^{l,n+1}_{ms} = (R^{l,n}_T)^T T^{l,n+1}_c, \\
 \ \ A^{l,n}_{p,c}=R^{l,n}_p A^{l,n}_p (R^{l,n}_p)^T, \ \ G_{p,c} = R^{l,n}_p G_p, \ \ p^{l,n+1}_{ms} = (R^{l,n}_p)^T p^{l,n+1}_c,
\end{split} \nonumber
\end{equation}}
and
\begin{equation}
\begin{split}
R^{l,n}_T=(\psi^1_{T,1},...,\psi^1_{T,N_c},...,\psi^{N_c}_{T,1},...,\psi^{N_c}_{T,N_c},\psi^1_{add},..,\psi^{N_p}_{add}, \vartheta^1_{T,l},...,\vartheta^{N_c}_{T,l})^T, \\
R^n_p=R^{L,n}_p=(\psi^1_{p,1},...,\psi^1_{p,M_{N_c}},...,\psi^{N_c}_{p,1},...,\psi^{N_c}_{p,M_{N_c}},\vartheta^1_{p,l},...,\vartheta^{N_c}_{p,l})^T,
\end{split} \nonumber
\end{equation}
with $l=1,...,L$, $R^n_T=R^{L,n}_T$, $R^n_p=R^{L,n}_p$.

\section{Numerical results}

\begin{figure}[h!]
\begin{center}
\includegraphics[width=0.4\linewidth]{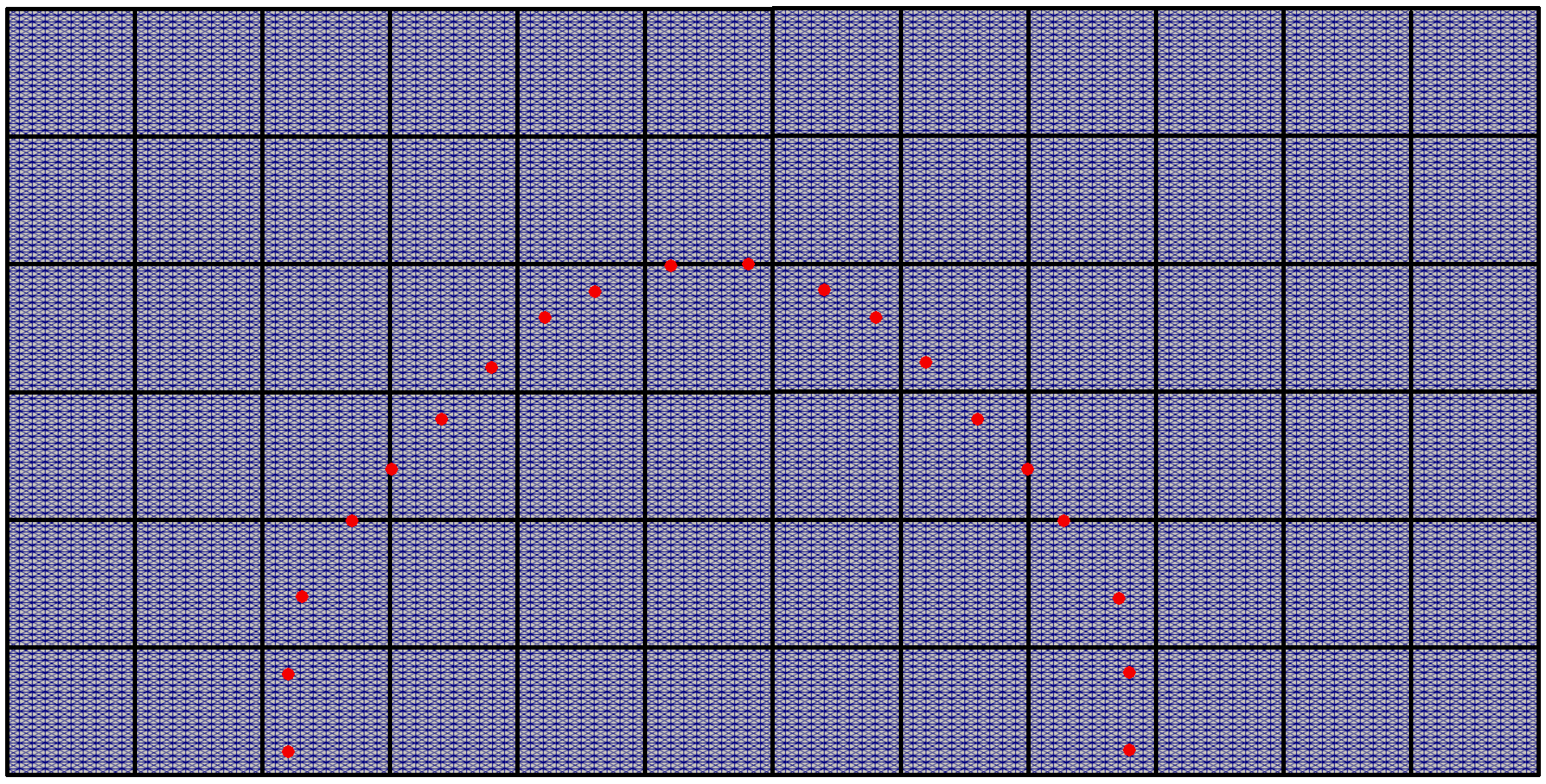}
\includegraphics[width=0.4\linewidth]{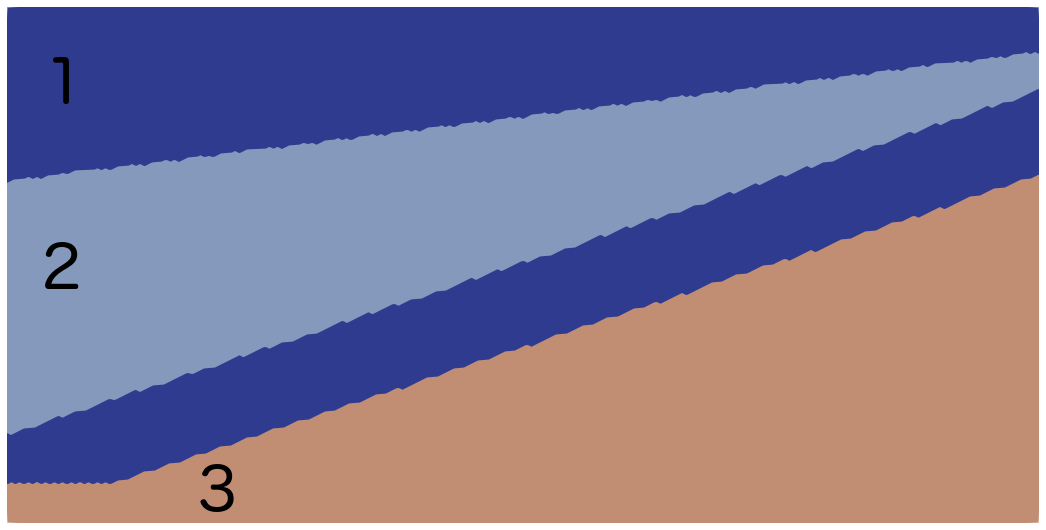}
\end{center}
\caption{The computational grid with 20 freezing pipes and heterogeneity. Left: coarse grid and fine grid. Right: heterogeneous properties.}
\label{domain}
\end{figure}

In this paper, we consider problems in two-dimensional formulation in heterogeneous domain with high contrast. We consider modeling of heat and mass transfer in computational domain $\Omega = [12\times 6][m]$. For modeling, we take $20$ freezing pipes in computational domain. Our heterogeneity is represented by layers and simulates a layered reservoir. We construct structured triangular fine grid $[120 \times 120]$ with 14641 vertices and 14400 cells. For multiscale solver, we generate coarse grid $24 \times 12$ with 325 vertices and 288 cells.   \revs{Computational domain, computational grid and schematic heterogeneity are presented in fig. \ref{domain}.} 

	We perform modeling of ground freezing process with $t_{max}=25$ days and $80$ times steps. In our simulation we set Dirichlet boundary condition on the freezing pipes $T_p=-30^\circ C$ and on the other boundaries we set $-k\frac{\partial T}{\partial n}=0$. As initial condition  for temperature we take $T_0=2^\circ C$. For pressure we take different boundary conditions:
	
\begin{itemize}

\item \textit{Test 1.} Flow from right to left, at the left boundary we set $p=1$, on the right boundary -- $p=0$, 

\item \textit{Test 2.} Flow from right to left, at the top boundary we set $p=1$, on the bottom boundary -- $p=0$.
 
\end{itemize} 
 On the other boundaries for pressure, we apply $-\lambda_\varepsilon  \frac{\partial p}{\partial n}= 0$.
 
We take the next values of coefficients:
\begin{equation}  \label{vals}
\begin{split}
\lambda ^+ _1 = 1.37 [W/(m\cdot K)], \ \lambda ^-_1 = 1.72 [W/(m\cdot K)], \ c^+_1 \rho^+_1 = 2.397\cdot 10^6 [J/(m^3 \cdot K],  \\ 
\ c^-_1 \rho^-_1 = 1.886 \cdot 10^6 [J/(m^3 \cdot K], \ \rho^+_1 L_1 = 75.33 \cdot 10^6, \ \frac{\rho _{w,1} \kappa _1}{\mu _1} = 1.0 \cdot 10 ^ {-13},  \\
\lambda ^+ _2 = 2.67 [W/(m\cdot K)], \ \lambda ^-_2 = 3.37 [W/(m\cdot K)], \ c^+_2 \rho^+_2 = 2.13 \cdot 10^6 [J/(m^3 \cdot K],  \\ 
\ c^-_2 \rho^-_2 = 2.09 \cdot 10^6 [J/(m^3 \cdot K], \ \rho^+_2 L_2 = 64.769 \cdot 10^6, \ \frac{\rho _{w,2} \kappa _2}{\mu _2} = 10.0 \cdot 10 ^ {-13},  \\
\lambda ^+ _3 = 1.4 [W/(m\cdot K)], \ \lambda ^-_3 = 1.56 [W/(m\cdot K)], \ c^+_3 \rho^+_3 = 2.96 \cdot 10^6 [J/(m^3 \cdot K],  \\ 
\ c^-_3 \rho^-_3 = 2.70 \cdot 10^6 [J/(m^3 \cdot K], \ \rho^+_2 L_2 = 130.544 \cdot 10^6, \ \frac{\rho _{w,2} \kappa _2}{\mu _2} = 5.0 \cdot 10 ^ {-13},  \\
\varepsilon = 10^{-3}, \ \triangle = 0.5, \ Q=0, \ F=0. \ \ \ \ \ \ \ \ \ \ \ \ \ \ \ \ \ \ \ \ \ \ \ \ \ \ \ \ \  \ \  \nonumber
\end{split}
\end{equation} 

\begin{figure}[h!]
\begin{center}
\includegraphics[width=0.8\linewidth]{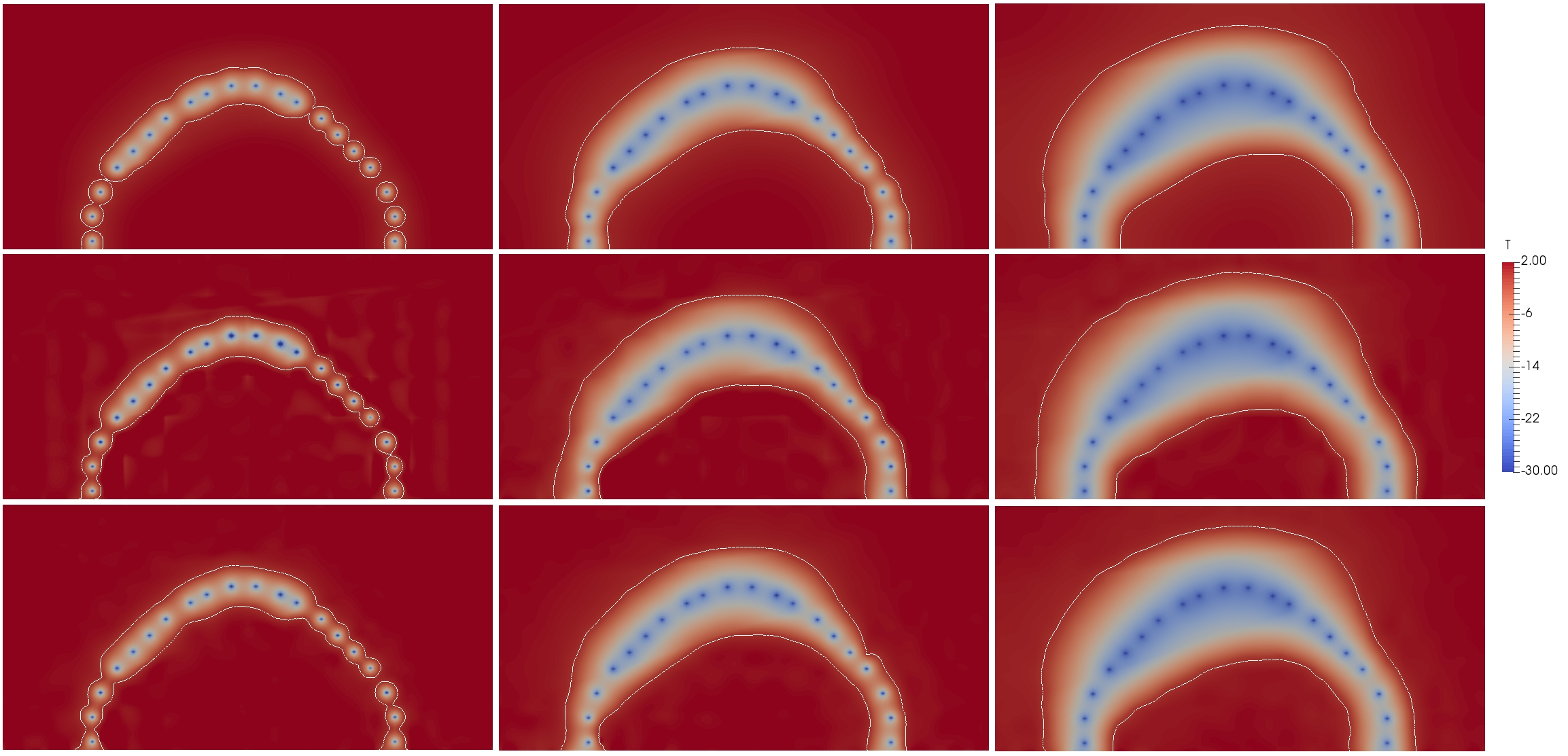}
\end{center}
\caption{Distribution of Temperature for \textit{Test 1} on 5, 30 and 80 time layers. First row: fine grid solution. 
Second row: multiscale solution using 4 offline basis functions.
Third row: multiscale solution using 4 offline basis functions and 1 online basis function.}
\label{result1temp}
\end{figure}
 
We will examine the multiscale and fine grid solutions for numerical results validation. To do this, we use relative errors in $L_2$ and $H_1$ norms:
\begin{equation}\label{erreq}
||e||_{L_2} = \sqrt{\frac{\int_\Omega (u_f - u_{ms})^2 dx}{\int_\Omega u^2_f dx}}, \ \ \  ||e||_{H_1} = \sqrt{\frac{\int_\Omega (\nabla (u_f - u_{ms}), \nabla (u_f - u_{ms})) dx}{\int_\Omega (\nabla u_f, \nabla u_f) dx}},
\end{equation}
where $u_f$ is a solution on a fine grid obtained by finite element method, $u_{ms}$ is a multiscale solution, instead of function $u$ we use temperature $T$ and pressure $p$.

\begin{figure}[h!]
\begin{center}
\includegraphics[width=0.8\linewidth]{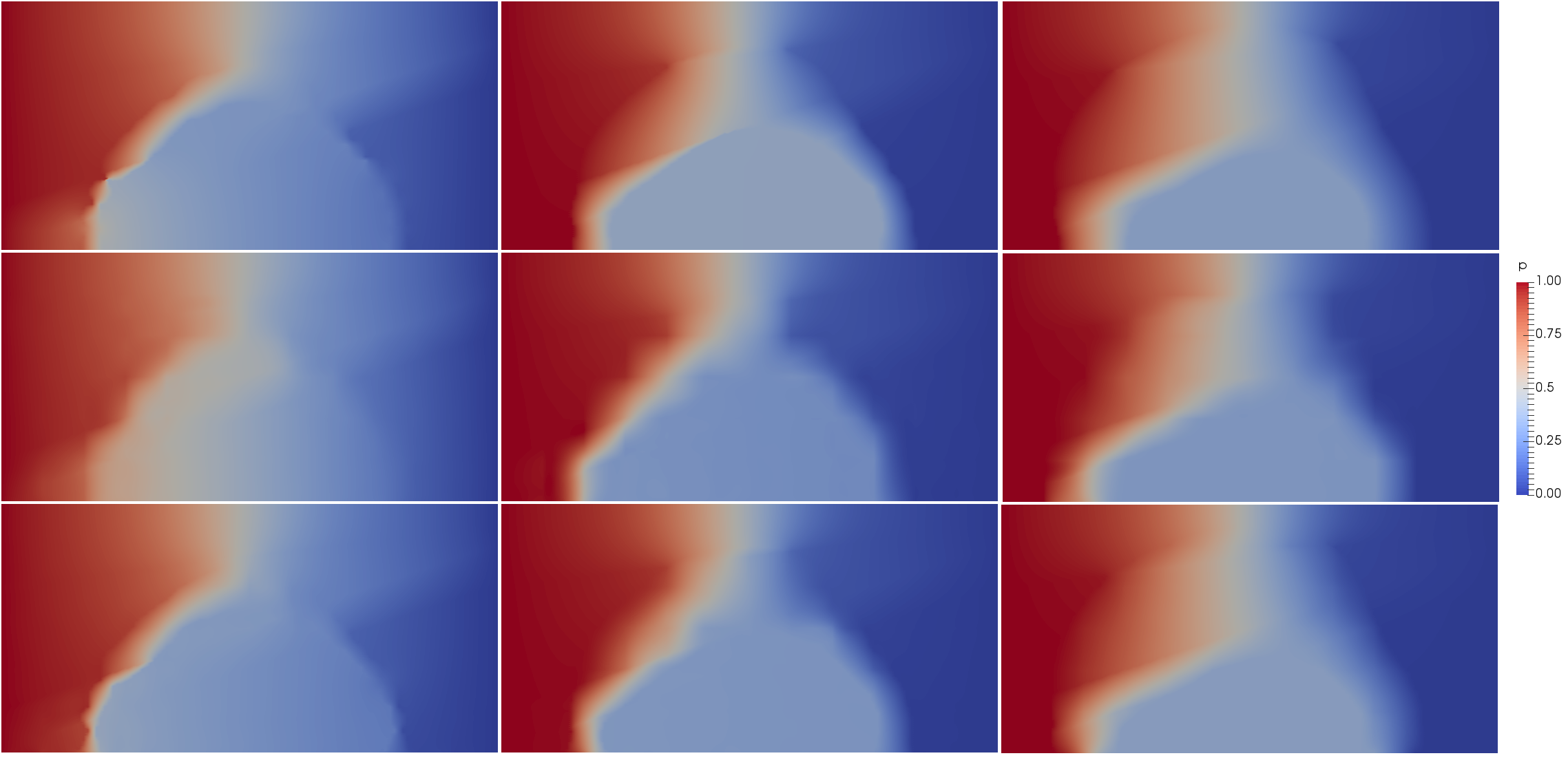}
\end{center}
\caption{Distribution of Pressure for \textit{Test 1} on 5, 30 and 80 time layers. First row: fine grid solution. 
Second row: multiscale solution using 4 offline basis functions.
Third row: multiscale solution using 4 offline basis functions and 1 online basis function.}
\label{result1press}
\end{figure}

\begin{table}[h!]
\begin{center}
\begin{tabular}{|c|c|cc|c|cc|c|cc|}
\hline
 Number of & \multicolumn{3}{|c|}{Offline basis} & \multicolumn{3}{|c|}{1 Online basis} & \multicolumn{3}{|c|}{2 Online basis} \\
\cline{2-10}
offline basis  
& \multirow{2}{*}{$DOF_c$ }  
& \multirow{2}{*}{$e^{L_2}_T$} & \multirow{2}{*}{$e^{H_1}_T$} & \multirow{2}{*}{$DOF_c$ }  & \multirow{2}{*}{$e^{L_2}_T$} & \multirow{2}{*}{$e^{H_1}_T$}  & \multirow{2}{*}{$DOF_c$ }  
& \multirow{2}{*}{$e^{L_2}_T$} & \multirow{2}{*}{$e^{H_1}_T$} \\
functions & & & & & & & & & \\
\hline
2 & 397 & 11.701 & 31.393 & 579 & 6.111 & 18.639 & 761 & 4.728 & 11.606 \\
4 & 761 & 8.222 & 21.188 & 943 & 5.484 & 13.826 & 1125 & 4.221 & 11.266 \\
6 & 1125 & 7.159 & 18.023 & 1307 & 4.893 & 12.438 & 1489 & 3.537 & 10.031 \\
8 & 1489 & 5.356 & 15.113 & 1671 & 3.757 & 10.174 & 1853 & 3.386 & 9.529 \\
\hline
\end{tabular}
\end{center}
\caption{Numerical results for \textit{Test 1}. Relative $L_2$ \revs{and $H_1$ ($\%$) errors for} Temperature with different number of offline multiscale basis functions using 1 and 2 online basis functions on the last time layer.}
\label{test1temp-err}
\end{table}

\begin{table}[h!]
\begin{center}
\begin{tabular}{|c|c|cc|c|cc|c|cc|}
\hline
 Number of & \multicolumn{3}{|c|}{Offline basis} & \multicolumn{3}{|c|}{1 Online basis} & \multicolumn{3}{|c|}{2 Online basis} \\
\cline{2-10}
offline basis  
& \multirow{2}{*}{$DOF_c$ }  
& \multirow{2}{*}{$e^{L_2}_p$} & \multirow{2}{*}{$e^{H_1}_p$} & \multirow{2}{*}{$DOF_c$ }  & \multirow{2}{*}{$e^{L_2}_p$} & \multirow{2}{*}{$e^{H_1}_p$}  & \multirow{2}{*}{$DOF_c$ }  
& \multirow{2}{*}{$e^{L_2}_p$} & \multirow{2}{*}{$e^{H_1}_p$} \\
functions & & & & & & & & & \\
\hline
2 & 397 & 15.817 & 66.365 & 579 & 1.831 & 28.653 & 761 & 1.149 & 21.409 \\
4 & 761 & 5.826 & 48.115 & 943 & 1.881 & 26.970 & 1125 & 1.556 & 24.346 \\
6 & 1125 & 4.603 & 44.995 & 1307 & 1.913 & 25.779 & 1489 & 1.519 & 22.885 \\
8 & 1489 & 3.799 & 40.171 & 1671 & 1.431 & 21.668 & 1853 & 1.221 & 19.637 \\
\hline
\end{tabular}
\end{center}
\caption{Numerical results for \textit{Test 1}. Relative $L_2$ \revs{and $H_1$ ($\%$) errors} for Pressure with different number of offline multiscale basis functions using 1 and 2 online basis functions on the last time layer.}
\label{test1press-err}
\end{table}

On several time layers, we display a temperature distribution \textit{Test 1} in Figure \ref{result1temp}, pressure distribution in time are shown in Figure \ref{result1temp}. In these figures in first row we present a fine grid solution, in second row we present a multiscale solution using 4 offline multiscale basis functions in each local domain and in third row we present a multiscale solution with 4 offline and 1 online multiscale basis functions. The additional basis functions \eqref{basisadd} is used in each computation in the set of offline basis functions. In Figure \ref{result1temp} the phase change boundary is drawn by a white line. One can see that our method describes the phase change boundary with high accuracy.  We should note, that we enrich multiscale space by online basis function on each 5-th time layer.

%\begin{figure}[h!]
%\begin{center}
%\includegraphics[width=0.485\linewidth]{figs/Test1L2ErrorTemp2.png}
%\ \
%\includegraphics[width=0.485\linewidth]{figs/Test1H1errorTemp2.png}
%\end{center}
%\caption{Numerical results for \textit{Test 1}. Error comparison for Temperature between offline approach using 2 offline basis and online approach using 1 and 2 online basis. Left: $L_2$ error. Right: $H_1$ error. }
%\label{Test1ErrorTemp}
%\end{figure}

%\begin{figure}[h!]
%\begin{center}
%\includegraphics[width=0.485\linewidth]{figs/Test1L2errorPress2.png}
%\ \
%\includegraphics[width=0.485\linewidth]{figs/Test1H1errorPress2.png}
%\end{center}
%\caption{Numerical results for \textit{Test 1}. Error comparison for Pressure between offline approach using 2 offline basis and online approach using 1 and 2 online basis. Left: $L_2$ error. Right: $H_1$ error. }
%\label{Test1ErrorTemp}
%\end{figure}

\begin{figure}[h!]
\begin{center}
\includegraphics[width=0.4\linewidth]{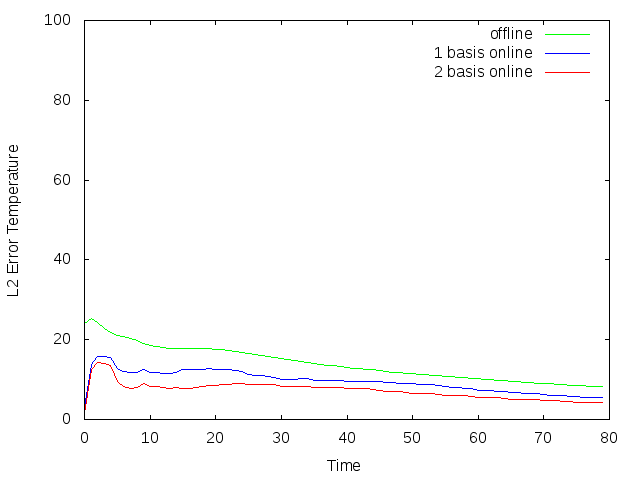}
\ \
\includegraphics[width=0.4\linewidth]{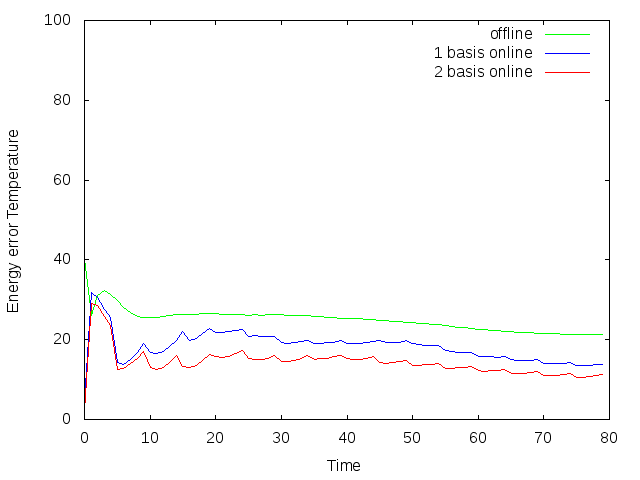}
\end{center}
\caption{Numerical results for \textit{Test 1}. Error comparison for Temperature between offline approach using 4 offline basis and online approach using 1 and 2 online basis. Left: $L_2$ error. Right: $H_1$ error. }
\label{Test1ErrorTemp}
\end{figure}

\begin{figure}[h!]
\begin{center}
\includegraphics[width=0.4\linewidth]{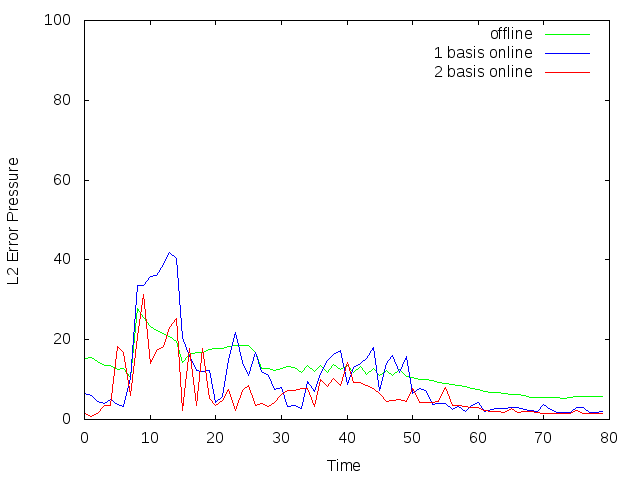}
\ \
\includegraphics[width=0.4\linewidth]{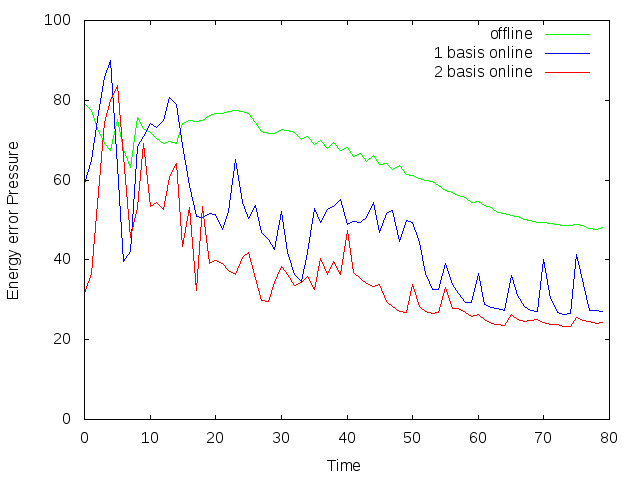}
\end{center}
\caption{Numerical results for \textit{Test 1}. Error comparison for Pressure between offline approach using 4 offline basis and online approach using 1 and 2 online basis. Left: $L_2$ error. Right: $H_1$ error. }
\label{Test1ErrorPress}
\end{figure}

\begin{figure}[h!]
\begin{center}
\includegraphics[width=0.4\linewidth]{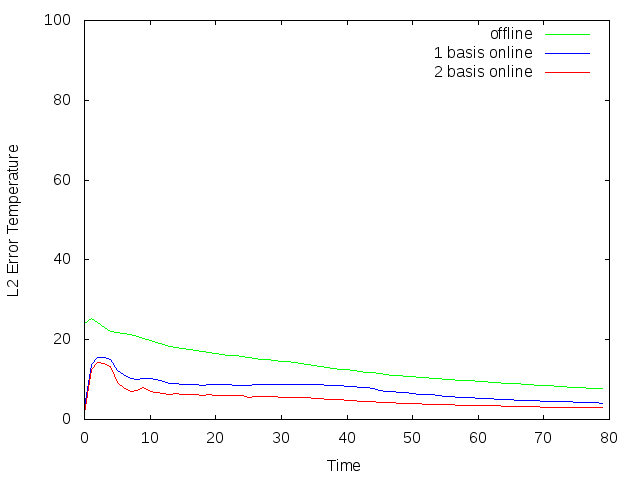}
\ \
\includegraphics[width=0.4\linewidth]{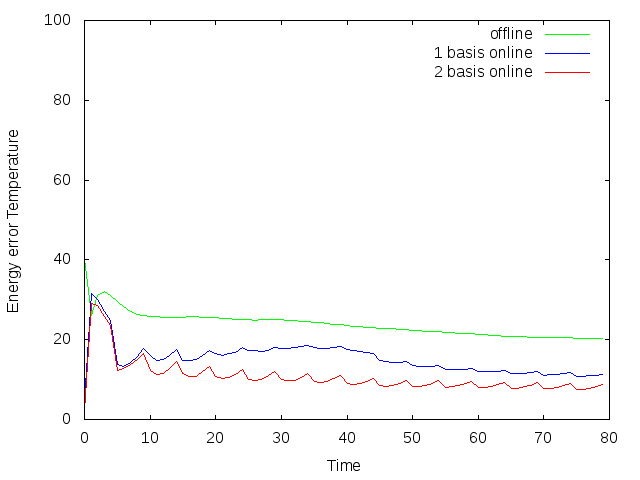}
\end{center}
\caption{Numerical results for \textit{Test 2}. Error comparison for Temperature between offline approach using 4 offline basis and online approach using 1 and 2 online basis. Left: $L_2$ error. Right: $H_1$ error. }
\label{Test2ErrorTemp4}
\end{figure}

We show the relative $L_2$ \revs{and $H_1$ errors for} temperature in Table \ref{test1temp-err} and for pressure in Table \ref{test1press-err}. In these tables, we show the comparison between simple GMsFEM and Online GMsFEM, using 1 and 2 online basis functions. From the tables, one can see, that we have a big error decay when we are using online basis functions. The proposed method has almost no difference in accuracy between using 1 and 2 online basis functions. We denote by $DOF_c$ a number of degree of freedom on the coarse grid. In Figures \ref{Test2ErrorTemp4} and \ref{Test1ErrorPress} we present $L_2$ and $H_1$ errors distribution in time for temperature and pressure respectively. The graphics present the errors for 4 offline multiscale basis functions with 0,1 or 2 online basis functions.  For online GMsFEM, we can observe the errors fall on each 5-th time layer. The error jumps are smaller when we use 2 online multiscale basis functions. From that, we can conclude that is better to use 2 online basis functions in each $\omega _i$.  

\begin{figure}[h!]
\begin{center}
\includegraphics[width=0.8\linewidth]{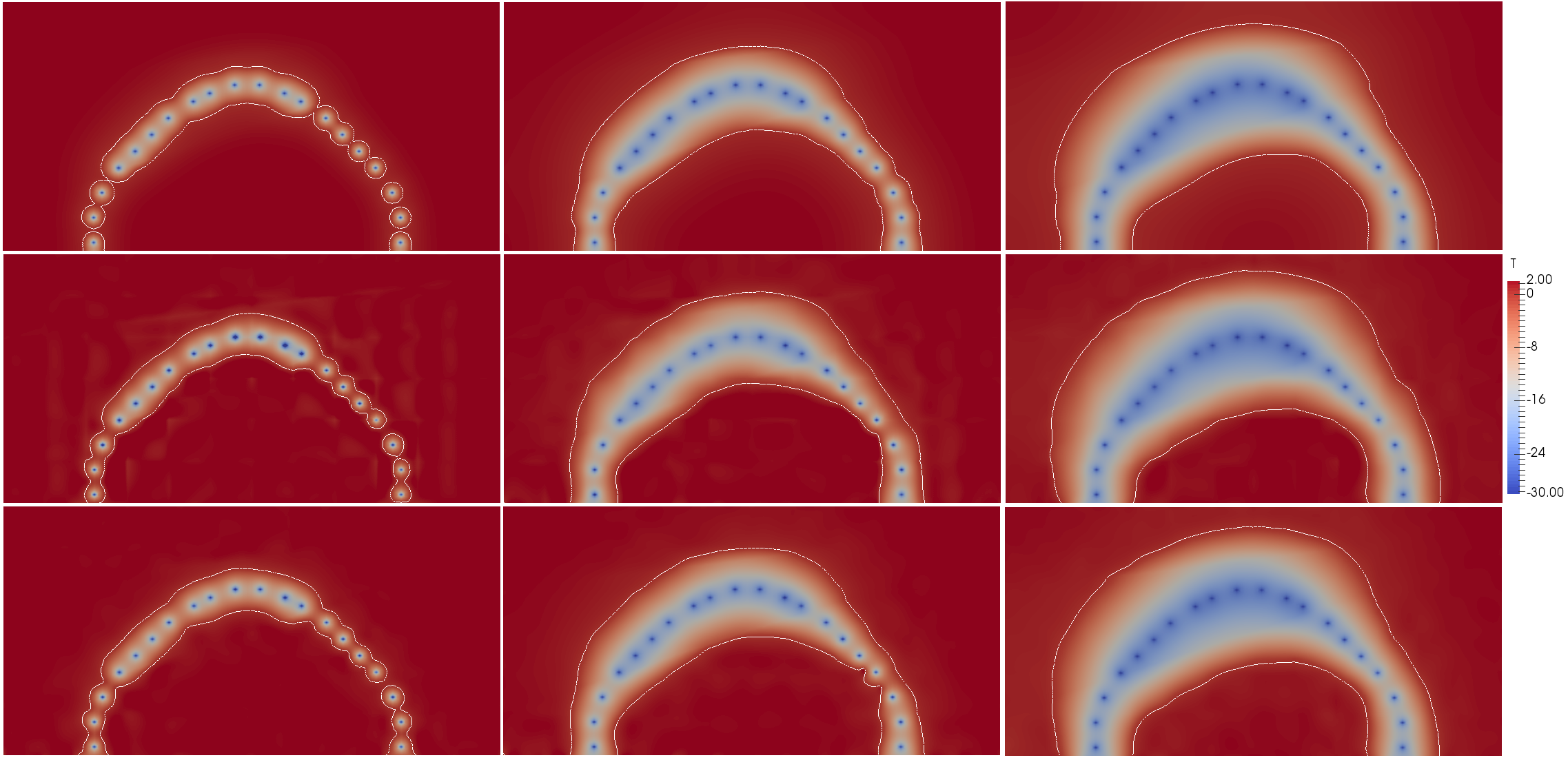}
\end{center}
\caption{Distribution of Temperature for \textit{Test 2} on 5, 30 and 80 time layers. First row: fine grid solution. 
Second row: multiscale solution using 4 offline basis functions.
Third row: multiscale solution using 4 offline basis functions and 1 online basis function.}
\label{result2temp}
\end{figure}

\begin{figure}[h!]
\begin{center}
\includegraphics[width=0.8\linewidth]{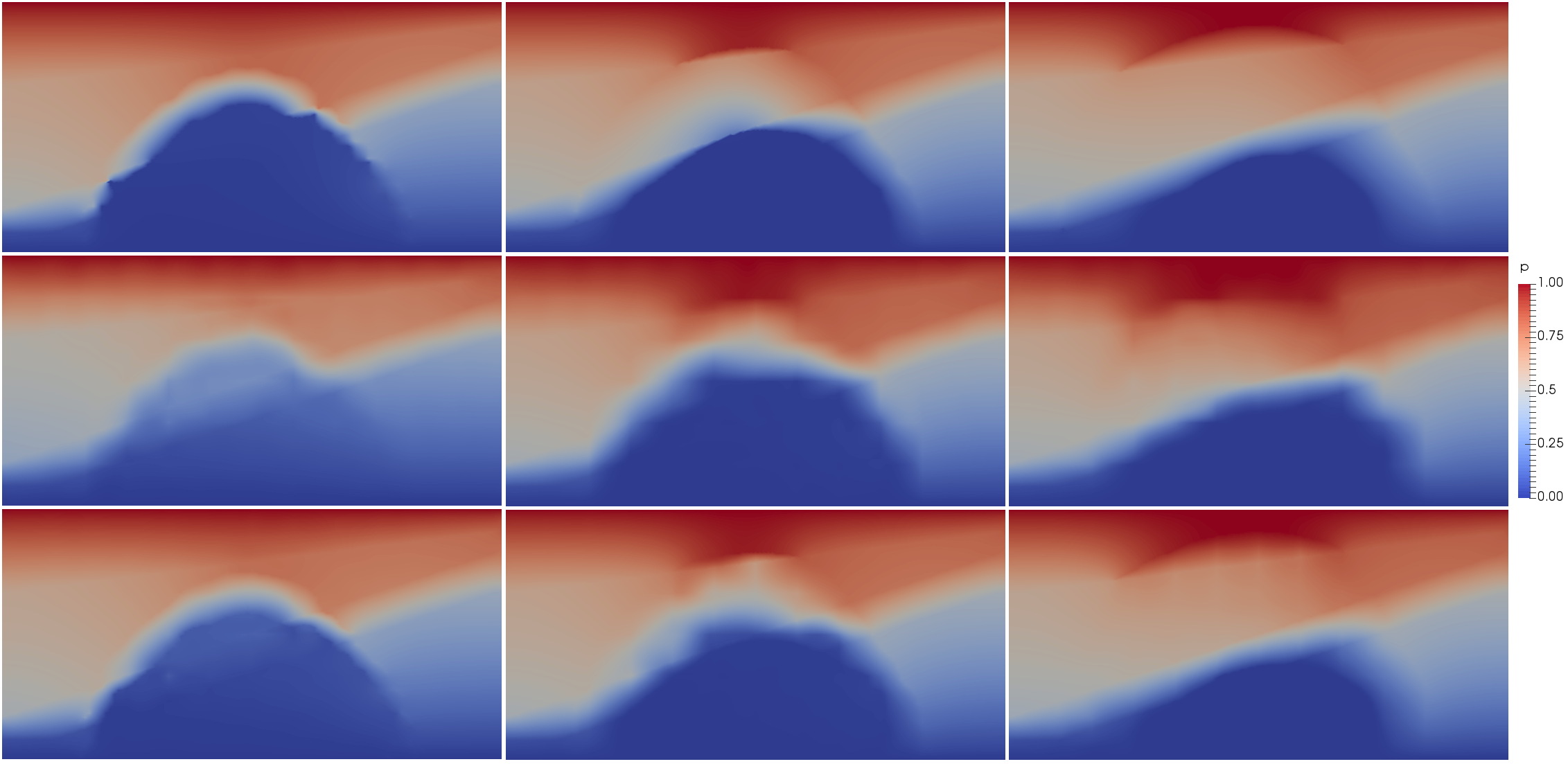}
\end{center}
\caption{Distribution of Pressure for \textit{Test 2} on 5, 30 and 80 time layers. First row: fine grid solution. 
Second row: multiscale solution using 4 offline basis functions.
Third row: multiscale solution using 4 offline basis functions and 1 online basis function.}
\label{result2press}
\end{figure}

For \textit{Test 2} we present numerical results in Figure \ref{result2temp} for temperature and in Figure \ref{result2press} for pressure. The multiscale solution is obtained by using 4 offline basis functions and additional basis. The fine grid solution is presented in the first row of the Figures, the offline solution is presented in the second row and in the third row we show the solution using 1 online multiscale basis functions. We can observe the good approximation of phase change boundary by Online GMsFEM. The change of flow direction in pressure didn't seem to affect the accuracy of the method. More accurate conclusions can be drawn by looking at the errors. 

\begin{table}[h!]
\begin{center}
\begin{tabular}{|c|c|cc|c|cc|c|cc|}
\hline
 Number of & \multicolumn{3}{|c|}{Offline basis} & \multicolumn{3}{|c|}{1 Online basis} & \multicolumn{3}{|c|}{2 Online basis} \\
\cline{2-10}
offline basis  
& \multirow{2}{*}{$DOF_c$ }  
& \multirow{2}{*}{$e^{L_2}_T$} & \multirow{2}{*}{$e^{H_1}_T$} & \multirow{2}{*}{$DOF_c$ }  & \multirow{2}{*}{$e^{L_2}_T$} & \multirow{2}{*}{$e^{H_1}_T$}  & \multirow{2}{*}{$DOF_c$ }  
& \multirow{2}{*}{$e^{L_2}_T$} & \multirow{2}{*}{$e^{H_1}_T$} \\
functions & & & & & & & & & \\
\hline
2 & 397 & 13.631 & 35.347 & 579 & 10.064  & 25.264 & 761 & 4.434 & 10.333 \\
4 & 761 & 7.709 & 20.231 & 943 & 4.129 & 11.337 & 1125 & 2.919 & 8.693 \\
6 & 1125 & 6.462 & 16.861 & 1307 & 3.819 & 10.510 & 1489 & 2.663 & 8.028 \\
8 & 1489 & 5.103 & 14.521 & 1671 & 3.021 & 8.857 & 1853 & 2.539 & 7.475 \\
\hline
\end{tabular}
\end{center}
\caption{Numerical results for \textit{Test 2}. Relative $L_2$ \revs{and $H_1$ ($\%$) errors} for Pressure with different number of offline multiscale basis functions using 1 and 2 online basis functions on the last time layer.}
\label{test2temp-err}
\end{table}

\begin{figure}[h!]
\begin{center}
\includegraphics[width=0.45\linewidth]{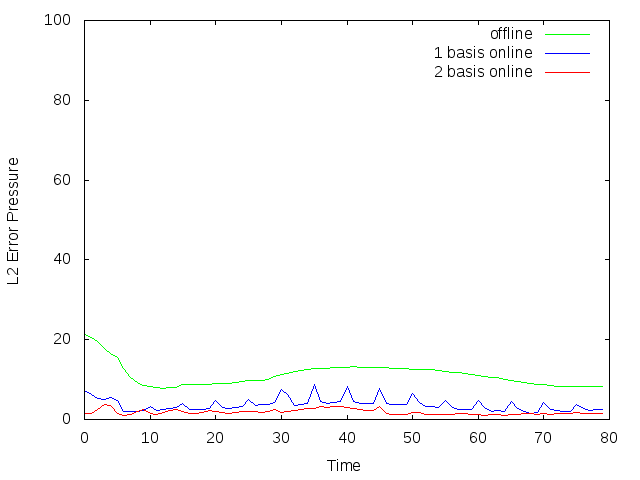}
\ \
\includegraphics[width=0.45\linewidth]{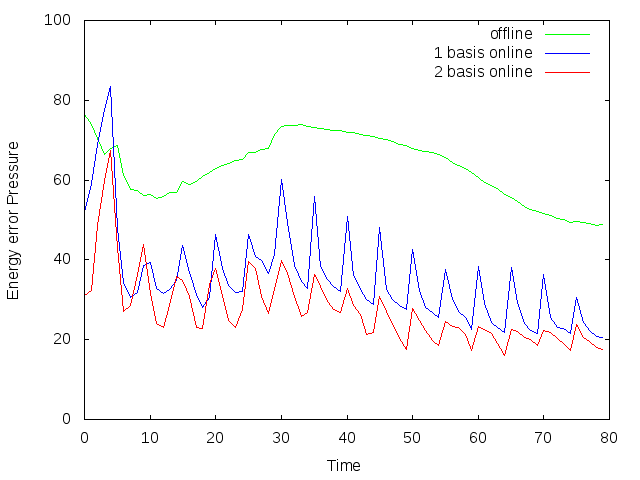}
\end{center}
\caption{Numerical results for \textit{Test 2}. Error comparison for Pressure between offline approach using 4 offline basis and online approach using 1 and 2 online basis. Left: $L_2$ error. Right: $H_1$ error. }
\label{Test2ErrorPress4}
\end{figure}

The relative $L_2$ \revs{and $H_1$ errors for} \textit{Test 2} are presented in Tables \ref{test2temp-err} and \ref{test2presserr}. The errors are shown for offline solution using 4 offline basis functions and for solutions with 1 and 2 online basis functions. The behavior of the errors is the same from \textit{Test 1}.  But here we have a little difference in behavior using online bases. Here we have a bigger difference in temperature error between the results with 1 and 2 online basis functions. One can see the error falling with adding second online basis more clearly. As in \textit{Test1}, the online approach allows taking less number of offline multiscale basis functions.  We also demonstrate the error distribution in time with 4 offline basis functions using 1 and 2 online bases in Figure \ref{Test2ErrorTemp4} for temperature and for pressure in Figure \ref{Test2ErrorPress4}. We also can see the error decay on each 5-th time layer. The jump of the error is smaller when we use 2 online basis functions. In this case, the adding of second online basis functions is also justified. 

\begin{table}[h!]
\begin{center}
\begin{tabular}{|c|c|cc|c|cc|c|cc|}
\hline
 Number of & \multicolumn{3}{|c|}{Offline basis} & \multicolumn{3}{|c|}{1 Online basis} & \multicolumn{3}{|c|}{2 Online basis} \\
\cline{2-10}
offline basis  
& \multirow{2}{*}{$DOF_c$ }  
& \multirow{2}{*}{$e^{L_2}_p$} & \multirow{2}{*}{$e^{H_1}_p$} & \multirow{2}{*}{$DOF_c$ }  & \multirow{2}{*}{$e^{L_2}_p$} & \multirow{2}{*}{$e^{H_1}_p$}  & \multirow{2}{*}{$DOF_c$ }  
& \multirow{2}{*}{$e^{L_2}_p$} & \multirow{2}{*}{$e^{H_1}_p$} \\
functions & & & & & & & & & \\
\hline
2 & 397 & 13.009 & 58.768 & 579 & 4.786 & 25.264 & 761 & 2.741 & 20.211 \\
4 & 761 & 8.341 & 48.899 & 943 & 2.491 & 20.659 & 1125 & 1.524 & 17.468 \\
6 & 1125 & 7.292 & 46.295 & 1307 & 1.801 & 17.433 & 1489 & 1.124 & 13.675 \\
8 & 1489 & 5.811 & 41.757 & 1671 & 1.371 & 14.712 & 1853 & 1.071 & 12.517 \\
\hline
\end{tabular}
\end{center}
\caption{Numerical results for \textit{Test 2}. Relative $L_2$ \revs{and $H_1$ ($\%$) errors for} Pressure with different number of offline multiscale basis functions using 1 and 2 online basis functions on the last time layer.}
\label{test2presserr}
\end{table}

We need to discuss obtained results. The Online GMsFEM showed a good results in both cases. \revs{ We should note that we did not show the results without an additional basis functions due to the fact that these results were obtained by us in our previous work \cite{vasilyeva2020multiscale}, which shows the effectiveness of the additional bases.} We need to decide about the number of online basis functions. And we need to mention In the $\textit{Test 1}$ we have a little error fall with adding second online basis, but in $\textit{Test 2}$ we have a bigger error falling for temperature. But in both cases the second online basis reduced the error jump in time. So, the question of adding a second basis is for the desired purposes. If for a certain study we want to know only the final state, then just enough only 1 online basis. If we need to know the solution at every moment of time, then we should use 2 online bases.

%\begin{figure}[h!]
%\begin{center}
%\includegraphics[width=0.485\linewidth]{figs/Test2L2ErrorTemp2.png}
%\ \
%\includegraphics[width=0.485\linewidth]{figs/Test2H1ErrorTemp2.png}
%\end{center}
%\caption{Numerical results for \textit{Test 2}. Error comparison for Temperature between offline approach using 2 offline basis and online approach using 1 and 2 online basis. Left: $L_2$ error. Right: $H_1$ error. }
%\label{Test2ErrorTemp2}
%\end{figure}

%\begin{figure}[h!]
%\begin{center}
%\includegraphics[width=0.485\linewidth]{figs/Test2L2ErrorPress2.png}
%\ \
%\includegraphics[width=0.485\linewidth]{figs/Test2H1ErrorPress2.png}
%\end{center}
%\caption{Numerical results for \textit{Test 2}. Error comparison for Pressure between offline approach using 2 offline basis and online approach using 1 and 2 online basis. Left: $L_2$ error. Right: $H_1$ error. }
%\label{Test2ErrorTemp2}
%\end{figure}

\section{Conclusion}

An online generalized multiscale finite element method for heat and mass problem in heterogeneous media with artificial ground freezing was proposed in this study. We considered two test cases with different boundary condition for pressure. We did a numerical experiment to demonstrate the method's correctness. We demonstrated multiscale solutions with various numbers of offline and online multiscale basis functions during the experiment. The experiment clearly showed that the addition of online basis function greatly improves the accuracy of GMsFEM. But whether to add a second online basis function already depends on the goals of the study. Nevertheless, in all the performed experiments, we obtained good error rates. We have significant decrease of the original system while not losing in terms of accuracy. An online generalized multiscale finite element method can be used for modeling heat and mass transfer problems.

\section{Acknowledgements}

This work is supported the grant of Russian Science Foundation N21-71-00061 and the Russian government project Science and Universities (project FSRG-2021-0015) aimed at supporting junior laboratories.

\bibliographystyle{unsrt}
\bibliography{lit}

\end{document}